# OPTIMAL SMOOTHING IN NONPARAMETRIC MIXED-EFFECT MODELS[1]

By Chong Gu and Ping Ma

*Purdue University and Harvard University*

Mixed-effect models are widely used for the analysis of correlated data such as longitudinal data and repeated measures. In this article, we study an approach to the nonparametric estimation of mixed-effect models. We consider models with parametric random effects and flexible fixed effects, and employ the penalized least squares method to estimate the models. The issue to be addressed is the selection of smoothing parameters through the generalized cross-validation method, which is shown to yield optimal smoothing for both real and latent random effects. Simulation studies are conducted to investigate the empirical performance of generalized cross-validation in the context. Real-data examples are presented to demonstrate the applications of the methodology.

**1. Introduction.** Mixed-effect models are widely used for the analysis of data with correlated errors. The linear mixed-effect models, also known as variance component models, are of the form

$$Y_i = \mathbf{x}_i^T \boldsymbol{\beta} + \mathbf{z}_i^T \mathbf{b} + \varepsilon_i, \qquad (1.1)$$

$i = 1, \ldots, n$, where $\mathbf{x}_i^T \boldsymbol{\beta}$ are the fixed effects, $\mathbf{z}_i^T \mathbf{b}$ are the random effects with $\mathbf{b} \sim N(\mathbf{0}, B)$, and $\varepsilon_i \sim N(0, \sigma^2)$ are independent of $\mathbf{b}$ and of each other; see, for example, [5] and [12]. The unknown parameters are $\boldsymbol{\beta}$, $B$ and $\sigma^2$, which are to be estimated from the data. Nonlinear and nonparametric generalizations of (1.1) can be found in, for example, [8, 11, 17].

In this article, we consider models of the form

$$Y_i = \eta(x_i) + \mathbf{z}_i^T \mathbf{b} + \varepsilon_i, \qquad (1.2)$$

Received May 2002; revised August 2004.
[1]Supported by NIH Grant R33HL68515.

where the regression function $\eta(x)$ is assumed to be a smooth function on a generic domain $\mathcal{X}$. The model terms $\eta(x)$ or $\eta(x) + \mathbf{z}^T\mathbf{b}$ will be estimated using the penalized (unweighted) least squares method through the minimization of

$$(1.3) \qquad \frac{1}{n}\sum_{i=1}^{n}(Y_i - \eta(x_i) - \mathbf{z}_i^T\mathbf{b})^2 + \frac{1}{n}\mathbf{b}^T\Sigma\mathbf{b} + \lambda J(\eta),$$

where the quadratic functional $J(\eta)$ quantifies the roughness of $\eta$ and the smoothing parameter $\lambda$ controls the trade-off between the goodness-of-fit and the smoothness of $\eta$; note that if one substitutes $\sigma^2 B^{-1}$ for $\Sigma$ in (1.3), then the first two terms are proportional to the minus log likelihood of $(\mathbf{Y}, \mathbf{b})$. We will treat $\Sigma$ as a tuning parameter like $\lambda$, however, and not be concerned with the estimation of $\sigma^2 B^{-1}$. Technically, (1.3) resembles a partial spline model, but with the partial terms $\mathbf{z}^T\mathbf{b}$ penalized.

Absent the random effects $\mathbf{z}^T\mathbf{b}$, penalized least squares regression has been studied extensively in the literature; see, for example, [16] and [2] for comprehensive treatments of the subject. The models of (1.2) were first considered by Wang [17], who used penalized marginal likelihood (of $\mathbf{Y}$) to estimate $\eta$. Smoothing parameter selection in penalized marginal likelihood estimation with correlated data was studied by Wang [18], who illustrated the middling performance of various versions of cross-validation, in contrast to the more reliable performance of the generalized maximum likelihood method of Wahba [15] derived under the Bayes model of smoothing splines. Under the Bayes model, $\eta$ itself is decomposed into fixed and random effects, with $\lambda J(\eta)$ acting as the minus log likelihood of the random effects; the generalized maximum likelihood method of Wahba [15] is essentially the popular restricted maximum likelihood method widely used for the estimation of variance component models.

The purpose of this article is to study the estimation of the model terms in (1.2) through the minimization of (1.3), with the smoothing parameter $\lambda$ and the correlation parameters $\Sigma$ selected by the standard generalized cross-validation method of Craven and Wahba [1], which was developed for independent data. In some applications, the random effects $\mathbf{z}^T\mathbf{b}$ are physically interpretable, or real, and in some others, $\mathbf{z}^T\mathbf{b}$ are merely a convenient device for the modeling of variance components, or latent; for the latter case, the estimation through (1.3) turns the variance components into "mean components." For both real and latent random effects, generalized cross-validation will be shown to yield optimal smoothing, through asymptotic analysis and numerical simulation. Real-data examples are also presented to illustrate the applications of the methodology.

The rest of this article is organized as follows. In Section 2 the problem is formulated and preliminary analysis is conducted. Examples are given in Section 3. Generalized cross-validation and its optimality are discussed in Section 4, followed by simulation studies in Section 5. Real-data examples are shown in Section 6. Proofs of the theorems and lemmas in Section 4 are collected in Section 7. A few remarks in Section 8 conclude the article.

**2. Penalized least squares estimation.** Consider the minimization of (1.3) for $\eta$ in a $q$-dimensional space span$\{\xi_1,\ldots,\xi_q\}$. Functions in the space can be expressed as

$$\eta(x) = \sum_{j=1}^{q} c_j \xi_j(x) = \boldsymbol{\xi}^T(x)\mathbf{c}. \tag{2.1}$$

Plugging (2.1) into (1.3), one minimizes

$$(\mathbf{Y} - R\mathbf{c} - Z\mathbf{b})^T(\mathbf{Y} - R\mathbf{c} - Z\mathbf{b}) + \mathbf{b}^T \Sigma \mathbf{b} + n\lambda \mathbf{c}^T Q \mathbf{c} \tag{2.2}$$

with respect to $\mathbf{c}$ and $\mathbf{b}$, where $\Sigma > 0$ is $p \times p$, $R$ is $n \times q$ with the $(i,j)$th entry $\xi_j(x_i)$, $Z = (\mathbf{z}_1,\ldots,\mathbf{z}_n)^T$ is $n \times p$ and $Q$ is $q \times q$ with the $(j,k)$th entry $J(\xi_j, \xi_k)$. Differentiating (2.2) with respect to $\mathbf{c}$ and $\mathbf{b}$ and setting the derivatives to 0, one has

$$\begin{pmatrix} R^T R + n\lambda Q & R^T Z \\ Z^T R & Z^T Z + \Sigma \end{pmatrix} \begin{pmatrix} \mathbf{c} \\ \mathbf{b} \end{pmatrix} = \begin{pmatrix} R^T \mathbf{Y} \\ Z^T \mathbf{Y} \end{pmatrix}. \tag{2.3}$$

Assume that the linear system is solvable, that is, the columns of $\binom{R^T}{Z^T}$ are in the column space of the left-hand side matrix. A solution of (2.3) is then given by

$$\begin{pmatrix} \hat{\mathbf{c}} \\ \hat{\mathbf{b}} \end{pmatrix} = \begin{pmatrix} R^T R + n\lambda Q & R^T Z \\ Z^T R & Z^T Z + \Sigma \end{pmatrix}^+ \begin{pmatrix} R^T \mathbf{Y} \\ Z^T \mathbf{Y} \end{pmatrix},$$

where $C^+$ denotes the Moore–Penrose inverse of $C$ satisfying $CC^+C = C$, $C^+CC^+ = C^+$, $(CC^+)^T = CC^+$ and $(C^+C)^T = C^+C$.

Write $D = Z^T Z + \Sigma$ and $E = (R^T R + n\lambda Q) - R^T Z D^{-1} Z^T R$. With (2.3) solvable, one has

$$\begin{pmatrix} R^T R + n\lambda Q & R^T Z \\ Z^T R & D \end{pmatrix} \begin{pmatrix} K \\ L \end{pmatrix} = \begin{pmatrix} R^T \\ Z^T \end{pmatrix}$$

for some $K$ and $L$, which, after some algebra, yields $EK(I - ZD^{-1}Z^T)^{-1} = R^T$, so the columns of $R^T$ are in the column space of $E$. It follows that $EE^+R^T = R^T$, and in turn

$$\begin{pmatrix} R^T R + n\lambda Q & R^T Z \\ Z^T R & Z^T Z + \Sigma \end{pmatrix}^+$$
$$= \begin{pmatrix} E^+ & -E^+ R^T Z D^{-1} \\ -D^{-1} Z^T R E^+ & D^{-1} + D^{-1} Z^T R E^+ R^T Z D^{-1} \end{pmatrix}.$$

It then follows that

$$\hat{\boldsymbol{\eta}} = R\hat{\mathbf{c}} = RE^+ R^T (I - ZD^{-1}Z^T)\mathbf{Y} = M\mathbf{Y}. \tag{2.4}$$

Similarly, one has

$$\hat{\mathbf{Y}} = R\hat{\mathbf{c}} + Z\hat{\mathbf{b}}$$
$$= \{(I - ZD^{-1}Z^T)RE^+ R^T (I - ZD^{-1}Z^T) + ZD^{-1}Z^T\}\mathbf{Y} = A(\lambda, \Sigma)\mathbf{Y},$$

where

$$A(\lambda, \Sigma) = (R, Z) \begin{pmatrix} R^T R + n\lambda Q & R^T Z \\ Z^T R & Z^T Z + \Sigma \end{pmatrix}^+ \begin{pmatrix} R^T \\ Z^T \end{pmatrix}$$
$$= (I - ZD^{-1}Z^T)RE^+ R^T (I - ZD^{-1}Z^T) + ZD^{-1}Z^T \tag{2.5}$$

is known as the smoothing matrix. Alternatively, for $\tilde{E} = R^T R + n\lambda Q$ and $\tilde{D} = D - Z^T R \tilde{E}^+ R^T Z$, one may write

$$\begin{pmatrix} R^T R + n\lambda Q & R^T Z \\ Z^T R & Z^T Z + \Sigma \end{pmatrix}^+$$
$$= \begin{pmatrix} \tilde{E}^+ + \tilde{E}^+ R^T Z \tilde{D}^{-1} Z^T R \tilde{E}^+ & -\tilde{E}^+ R^T Z \tilde{D}^{-1} \\ -\tilde{D}^{-1} Z^T R \tilde{E}^+ & \tilde{D}^{-1} \end{pmatrix},$$

yielding the expressions

$$M = \tilde{A}(\lambda) - \tilde{A}(\lambda)Z(Z^T (I - \tilde{A}(\lambda))Z + \Sigma)^{-1} Z^T (I - \tilde{A}(\lambda)), \tag{2.6}$$

where $\tilde{A}(\lambda) = R\tilde{E}^+ R^T$ is the smoothing matrix when the random effects are absent, and

$$A(\lambda, \Sigma) = \tilde{A}(\lambda) + (I - \tilde{A}(\lambda))Z(Z^T (I - \tilde{A}(\lambda))Z + \Sigma)^{-1} Z^T (I - \tilde{A}(\lambda)). \tag{2.7}$$

The eigenvalues of $A(\lambda, \Sigma)$ and $\tilde{A}(\lambda)$ are in the range $[0, 1]$.

With the standard formulation of penalized least squares regression, the minimization of (1.3) is performed in a so-called reproducing kernel Hilbert space $\mathcal{H} \subseteq \{\eta : J(\eta) < \infty\}$ in which $J(\eta)$ is a square seminorm, and the solution resides in the space $\mathcal{N}_J \oplus \text{span}\{R_J(x_i, \cdot), i = 1, \ldots, n\}$, where $\mathcal{N}_J = \{\eta : J(\eta) = 0\}$ is the null space of $J(\eta)$ and $R_J(\cdot, \cdot)$ is the so-called reproducing kernel in $\mathcal{H} \ominus \mathcal{N}_J$. The solution has an expression

$$\eta(x) = \sum_{i=1}^{m} d_\nu \phi_\nu(x) + \sum_{i=1}^{n} \tilde{c}_i R_J(x_i, x), \tag{2.8}$$

where $\{\phi_\nu\}_{\nu=1}^m$ is a basis of $\mathcal{N}_J$. It follows that $R = (S, \tilde{Q})$, where $S$ is $n \times m$ with the $(i, \nu)$th entry $\phi_\nu(x_i)$ and $\tilde{Q}$ is $n \times n$ with the $(i, j)$th entry $R_J(x_i, x_j)$. From the property of reproducing kernels, it can also be shown that $J(R_J(x_i, \cdot), R_J(x_j, \cdot)) = R_J(x_i, x_j)$, so $Q = \text{diag}(O, \tilde{Q})$. See, for example, [2] and [16]. The linear system (2.3) is thus solvable as long as $Z$ is of full column rank.

For fast computation, Kim and Gu [9] consider the space $\mathcal{N}_J \oplus \mathrm{span}\{R_J(z_j, \cdot), j = 1, \ldots, \tilde{q}\}$, where $\{z_j\}$ are a random subset of $\{x_i\}$. In that setting, $R = (S, \tilde{R})$, where $\tilde{R}$ is $n \times \tilde{q}$ with the $(i,j)$th entry $R_J(z_j, x_i)$, and $Q = \mathrm{diag}(O, \tilde{Q})$, where $\tilde{Q}$ is $\tilde{q} \times \tilde{q}$ with the $(j,k)$th entry $R_J(z_j, z_k)$. Since $J(\eta)$ is a square norm in $\mathrm{span}\{R_J(z_j, \cdot), j = 1, \ldots, \tilde{q}\}$, it can be shown that the columns of $\tilde{R}^T$ are in the column space of $\tilde{Q}$. It then follows that the linear system (2.3) is solvable when $Z$ is of full column rank.

The formulation of (2.1) and (2.2) also covers general penalized regression splines, so long as (2.3) is solvable. A sufficient condition is for both $R$ and $Z$ to be of full column rank.

**3. Examples.** A few examples are in order to illustrate the formulation of the problem and the potential applications of the method under study. The examples will be employed in the simulation study of Section 5 and the data analysis of Section 6.

EXAMPLE 3.1 (*Growth curves*). Consider the "growth" over time of a certain quantity associated with $p$ subjects,

$$Y_i = \eta(x_i) + b_{s_i} + \varepsilon_i,$$

where $Y_i$ is the $i$th observation taken at time $x_i \in [0, a]$ from subject $s_i \in \{1, \ldots, p\}$, and $b_s \sim N(0, \sigma_s^2)$ are the subject random effects, independent of the measurement error $\varepsilon_i$ and of each other. In this setting, $B = \sigma_s^2 I$, so the $p \times p$ matrix $\Sigma$ is diagonal with only one tunable parameter. The random effects $b_s$ are real.

Taking $J(\eta) = \int_0^a (d^2\eta/dx^2)^2 \, dx$, one has the cubic smoothing spline, with the $\phi_\nu$ and $R_J$ functions in (2.8) given by

$$\phi_1(x) = 1, \qquad \phi_2(x) = x, \qquad R_J(x_1, x_2) = \int_0^a (x_1 - u)_+ (x_2 - u)_+ \, du,$$

where $(\cdot)_+ = \max(\cdot, 0)$. See, for example, [2], Section 2.3.1. The null space model has the expression $\eta(x) = \beta_0 + \beta_1 x$.

Taking $J(\eta) = \int_0^a (L_\theta \eta)^2 h_\theta \, dx$, where $L_\theta = (d/dx)(d/dx + \theta)$ and $h_\theta = e^{3\theta x}$ for some $\theta > 0$, one has a (negative) exponential spline. The null space model has the expression $\eta(x) = \beta_0 + \beta_1 e^{-\theta x}$. Transforming $x$ by $\tilde{x} = (1 - e^{-\theta x})/\theta$, it can be shown that

$$\int_0^a (L_\theta \eta)^2 h_\theta \, dx = \int_0^{\tilde{a}} (d^2\eta/d\tilde{x}^2)^2 \, d\tilde{x},$$

where $\tilde{a} = (1 - e^{-\theta a})/\theta$, so one has a cubic spline in $\tilde{x}$. See, for example, [2], Example 4.7, Section 4.3.4. Note that the exponential spline reduces to the cubic spline in $x$ when $\theta = 0$.

Suppose $Y$ is the logarithm of the measurement $\tilde{Y}$ satisfying a log-normal distribution with $\mu = \eta(x) + b_s$ and $\sigma^2$ a constant; the mean of $\tilde{Y}$ is known to be $\exp(\mu + \sigma^2/2)$. The null space model of the cubic spline characterizes an exponential growth curve for $\tilde{Y}$, and the null space model of the exponential spline corresponds to a Gompertz growth curve for $\tilde{Y}$. The splines allow departures from these parametric growth curves.

EXAMPLE 3.2 (*Growth under treatment*). Consider the setting of Example 3.1, but with the $p$ subjects divided into $t$ treatment groups. The fixed effect becomes $\eta(x, \tau)$, where $\tau \in \{1, \ldots, t\}$ denotes the treatment level. For the identifiability of $\eta(x, \tau)$ and $b_s$, one needs more than one subject per treatment level. One may decompose

$$\eta(x, \tau) = \eta_\varnothing + \eta_1(x) + \eta_2(\tau) + \eta_{1,2}(x, \tau),$$

where $\eta_\varnothing$ is a constant, $\eta_1(x)$ is a function of $x$ satisfying $\eta_1(0) = 0$, $\eta_2(\tau)$ is a function of $\tau$ satisfying $\sum_{\tau=1}^t \eta_2(\tau) = 0$, and $\eta_{1,2}(x, \tau)$ satisfies $\eta_{1,2}(0, \tau) = 0$, $\forall \tau$, and $\sum_{\tau=1}^t \eta_{1,2}(x, \tau) = 0$, $\forall x$. The term $\eta_\varnothing + \eta_1(x)$ is the "average growth" and the term $\eta_2(\tau) + \eta_{1,2}(x, \tau)$ is the "contrast growth."

For flexible models one may use

$$J(\eta) = \theta_1^{-1} \int_0^a (d^2\eta_1/dx^2)^2 \, dx + \theta_{1,2}^{-1} \int_0^a \sum_{\tau=1}^t (d^2\eta_{1,2}/dx^2)^2 \, dx,$$

which has a null space $\mathcal{N}_J$ of dimension $2t$. A set of $\phi_\nu$ is given by

$$\{1, x, I_{[\tau=j]} - 1/t, (I_{[\tau=j]} - 1/t)x, j = 1, \ldots, t-1\},$$

and the function $R_J$ is given by

$$R_J(x_1, \tau_1; x_2, \tau_2) = \theta_1 \int_0^a (x_1 - u)_+ (x_2 - u)_+ \, du$$
$$+ \theta_{1,2}(I_{[\tau_1=\tau_2]} - 1/t) \int_0^a (x_1 - u)_+ (x_2 - u)_+ \, du.$$

See, for example, [2], Section 2.4.4, Problem 2.14(c). To force an additive model $\eta(x, \tau) = \eta_\varnothing + \eta_1(x) + \eta_2(\tau)$, which yields parallel growth curves at different treatment levels, one may set $\theta_{1,2} = 0$ and remove $(I_{[\tau=j]} - 1/t)x$ from the list of $\phi_\nu$. One may also choose to transform $x$ through $\tilde{x} = (1 - e^{-\theta x})/\theta$ and fit models on the $\tilde{x}$ scale.

EXAMPLE 3.3 (*Clustered observations*). Consider observations from $p$ clusters, such as in multicenter studies, $Y_i = \eta(x_i) + \tilde{\varepsilon}_i$, where $Y_i$ is taken from cluster $c_i$ with covariate $x_i$. Observations from different clusters are independent, while observations from the same cluster may be correlated to various degrees. The intracluster correlation may be modeled via $\tilde{\varepsilon}_i =$ $b_{c_i} + \varepsilon_i$, where $\mathbf{b} \sim N(0, B)$, with $B = \text{diag}(\sigma_1^2, \ldots, \sigma_p^2)$, and $\boldsymbol{\varepsilon} \sim N(0, \sigma^2 I)$, independent of each other; the intracluster correlation in cluster $c_i$ is given by $\sigma_i^2/(\sigma^2 + \sigma_i^2)$. In this setting, the $p \times p$ matrix $\Sigma$ involves $p$ tunable parameters on the diagonal. The random effects $b_c$ are latent.

Note that the covariate $x$ is generic, which can be univariate as in Example 3.1, or multivariate as in Example 3.2.

**4. Optimality of generalized cross-validation.** For the selection of the smoothing parameter $\lambda$ (and others such as the $\theta$ in Example 3.1 and the $\theta_1$ and $\theta_{1,2}$ in Example 3.2, if present) and the correlation parameters $\Sigma$, we propose to minimize the generalized cross-validation score

$$
(4.1) \qquad V(\lambda, \Sigma) = \frac{n^{-1} \mathbf{Y}^T (I - A(\lambda, \Sigma))^2 \mathbf{Y}}{\{n^{-1} \text{tr}(I - A(\lambda, \Sigma))\}^2};
$$

$\Sigma$ may involve less than $p(p+1)/2$ tunable parameters. It will be shown in this section that the minimizers of $V(\lambda, \Sigma)$ yield optimal smoothing asymptotically, in the sense to be specified. Numerical verifications of the asymptotic analysis will be presented in the next section. Generalized cross-validation was proposed by Craven and Wahba [1] for independent data, with the asymptotic optimality established by Li [10] in that setting; see also [13].

First consider the mean square error at the data points,

$$
(4.2) \qquad L_1(\lambda, \Sigma) = \frac{1}{n} \sum_{i=1}^{n} (\hat{Y}_i - \eta(x_i) - \mathbf{z}_i^T \mathbf{b})^2,
$$

which is a natural loss when the random effects $\mathbf{z}^T \mathbf{b}$ are real. Simple algebra yields

$$
\begin{aligned}
L_1(\lambda, \Sigma) &= \frac{1}{n}(A\mathbf{Y} - \boldsymbol{\eta} - Z\mathbf{b})^T (A\mathbf{Y} - \boldsymbol{\eta} - Z\mathbf{b}) \\
&= \frac{1}{n}(\boldsymbol{\eta} + Z\mathbf{b})^T (I - A)^2 (\boldsymbol{\eta} + Z\mathbf{b}) \\
&\quad - \frac{2}{n}(\boldsymbol{\eta} + Z\mathbf{b})^T (I - A) A \boldsymbol{\varepsilon} + \frac{1}{n} \boldsymbol{\varepsilon}^T A^2 \boldsymbol{\varepsilon},
\end{aligned}
$$

where $\boldsymbol{\eta} = (\eta(x_1), \ldots, \eta(x_n))^T$, $\mathbf{Y} = \boldsymbol{\eta} + Z\mathbf{b} + \boldsymbol{\varepsilon}$ and the arguments $(\lambda, \Sigma)$ are dropped from the notation of the smoothing matrix $A$. Taking expectation with respect to $\mathbf{b}$ and $\boldsymbol{\varepsilon}$, the risk is seen to be

$$
\begin{aligned}
(4.3) \quad R_1(\lambda, \Sigma) &= E[L_1(\lambda, \Sigma)] \\
&= \frac{1}{n} \boldsymbol{\eta}^T (I - A)^2 \boldsymbol{\eta} + \frac{1}{n} \text{tr}((I - A)^2 Z B Z^T) + \frac{\sigma^2}{n} \text{tr} A^2.
\end{aligned}
$$

Now define

(4.4) $$U(\lambda, \Sigma) = \frac{1}{n}\mathbf{Y}^T(I-A)^2\mathbf{Y} + \frac{2}{n}\sigma^2 \operatorname{tr} A.$$

It follows that

(4.5) $$\begin{aligned} U(\lambda, \Sigma) - L_1(\lambda, \Sigma) - \frac{1}{n}\varepsilon^T\varepsilon \\ = \frac{2}{n}(\boldsymbol{\eta} + Z\mathbf{b})^T(I-A)\varepsilon - \frac{2}{n}(\varepsilon^T A\varepsilon - \sigma^2 \operatorname{tr} A). \end{aligned}$$

We shall establish the optimality of $U(\lambda, \Sigma)$ under the following conditions.

CONDITION C.1. *The eigenvalues of $\Sigma(Z^T(I - \tilde{A}(\lambda))Z + \Sigma)^{-1}\Sigma$ are bounded from above.*

Condition C.1 holds for $\Sigma$ with eigenvalues bounded from above, and for $\Sigma$ of magnitude up to the order of $O(\sqrt{n})$ when the magnitude of $Z^T(I - \tilde{A}(\lambda))Z$ grows at the rate of $O(n)$.

CONDITION C.2. *As $n \to \infty$, $nR_1(\lambda, \Sigma) \to \infty$.*

The condition simply concedes that the parametric rate of $O(n^{-1})$ is not achievable. In the absence of random effects, for $\eta$ satisfying $J(\eta) < \infty$ or more stringent smoothness conditions, it typically holds that $n^{-1}\boldsymbol{\eta}^T(I - \tilde{A}(\lambda))^2\boldsymbol{\eta} = O(\lambda^s)$ for some $s \in [1,2]$, and $\operatorname{tr}\tilde{A}^2(\lambda) \asymp \lambda^{-1/r}$ as $\lambda \to 0$ and $n\lambda^{1/r} \to \infty$ for some $r > 1$, at least for univariate smoothing splines; see, for example, [1, 15] and [2], Section 4.2.3. For the cubic splines of Example 3.1, $r = 4$.

LEMMA 4.1. *Under Condition C.1, if $n^{-1}\boldsymbol{\eta}^T(I - \tilde{A}(\lambda))^2\boldsymbol{\eta} = O(\lambda^s)$ and $\operatorname{tr}\tilde{A}^2(\lambda) = O(\lambda^{-1/r})$ as $\lambda \to 0$ and $n\lambda^{1/r} \to \infty$, then $R_1(\lambda, \Sigma) = O(\lambda^s + n^{-1}\lambda^{-1/r} + n^{-1}p)$.*

See Section 7 for the proof of the lemma. For fixed $p$, the random effects add little to the equation, and Condition C.2 is satisfied for $\lambda \to 0$, $n\lambda^{1/r} \to \infty$ and $\Sigma$ of magnitude up to order $O(\sqrt{n})$; the optimal $\lambda \asymp n^{-r/(sr+1)}$ is well within the domain. In fact, the restriction on $\Sigma$ is not really necessary for Condition C.2 but to assure that $R_1 \to 0$. When $p$ grows with $n$, Condition C.2 clearly holds, though one may need to scale back the domain of $\Sigma$ for $R_1 = o(1)$ to remain true.

THEOREM 4.1. *Under Conditions C.1 and C.2, as $n \to \infty$, one has*
$$U(\lambda, \Sigma) - L_1(\lambda, \Sigma) - \frac{1}{n}\varepsilon^T\varepsilon = o_p(L_1(\lambda, \Sigma)).$$

The proof of the theorem is given in Section 7. When the conditions of the theorem hold in a neighborhood of the optimal $(\lambda, \Sigma)$, the minimizer of $U(\lambda, \Sigma)$ will deliver nearly the minimum loss.

The use of $U(\lambda, \Sigma)$ requires knowledge of $\sigma^2$, which usually is not available in practice. With an extra condition, the result also holds for $V(\lambda, \Sigma)$.

CONDITION C.3. As $n \to \infty$, $\{n^{-1} \operatorname{tr} A(\lambda, \Sigma)\}^2 / \{n^{-1} \operatorname{tr} A^2(\lambda, \Sigma)\} \to 0$.

In the absence of random effects, Condition C.3 generally holds in most settings of interest. In fact, it typically holds that $\operatorname{tr} \tilde{A}(\lambda) \asymp \lambda^{-1/r}$ as $\lambda \to 0$ and $n\lambda^{1/r} \to \infty$, of the same order as $\operatorname{tr} \tilde{A}^2(\lambda)$. See, for example, [1, 10, 15] and [2], Section 4.2.3.

LEMMA 4.2. *If* $\operatorname{tr} \tilde{A}(\lambda) = O(\lambda^{-1/r})$ *and* $\operatorname{tr} \tilde{A}^2(\lambda) \asymp \lambda^{-1/r}$ *as* $\lambda \to 0$ *and* $n\lambda^{1/r} \to \infty$, *then Condition* C.3 *holds for $p$ up to order* $O(\sqrt{n})$.

The proof is to be found in Section 7.

THEOREM 4.2. *Under Conditions* C.1, C.2 *and* C.3, *as $n \to \infty$, one has*

$$V(\lambda, \Sigma) - L_1(\lambda, \Sigma) - \frac{1}{n}\varepsilon^T \varepsilon = o_p(L_1(\lambda, \Sigma)).$$

PROOF. Given Theorem 4.1, the proof follows that of Theorem 3.3 in [2], page 66. □

We now turn to the case with latent random effects $\mathbf{z}^T \mathbf{b}$, for which the loss $L_1(\lambda, \Sigma)$ of (4.2) may not make much practical sense. Write $P_Z = Z(Z^T Z)^+ Z^T$ and $P_Z^\perp = I - P_Z$. We consider the estimation of $P_Z^\perp \boldsymbol{\eta}$ by $P_Z^\perp \hat{\boldsymbol{\eta}}$, where $\hat{\boldsymbol{\eta}}$ is given in (2.4); the projection ensures the identifiability of the target function. Accounting for the error covariance $\sigma^2 I + ZBZ^T$, one may assess the estimation precision via the loss

$$\tilde{L}_2(\lambda, \Sigma) = \frac{1}{n}(\hat{\boldsymbol{\eta}} - \boldsymbol{\eta})^T P_Z^\perp (\sigma^2 I + ZBZ^T)^{-1} P_Z^\perp (\hat{\boldsymbol{\eta}} - \boldsymbol{\eta}).$$

Since $(\sigma^2 I + ZBZ^T)^{-1} = \sigma^{-2}(I - ZD_0^{-1}Z^T)$, where $D_0 = Z^T Z + \sigma^2 B^{-1}$, one may use

(4.6) $\qquad L_2(\lambda, \Sigma) = \sigma^2 \tilde{L}_2(\lambda, \Sigma) = \frac{1}{n}(\hat{\boldsymbol{\eta}} - \boldsymbol{\eta})^T P_Z^\perp (\hat{\boldsymbol{\eta}} - \boldsymbol{\eta}),$

which is independent of $B$. Write $Q_Z = ZD^{-1}Z^T$ and recall $M = RE^+ R^T (I - Q_Z)$ from (2.4). Plugging $\hat{\boldsymbol{\eta}} = M(\boldsymbol{\eta} + Z\mathbf{b} + \varepsilon)$ into (4.6) and taking expectation, one has the risk

$$R_2(\lambda, \Sigma) = E[L_2(\lambda, \Sigma)]$$

$$\text{(4.7)} \qquad = \frac{1}{n}\{\boldsymbol{\eta}^T(I-M)^T P_Z^\perp (I-M)\boldsymbol{\eta}$$
$$+ \operatorname{tr}(M^T P_Z^\perp M Z B Z^T) + \sigma^2 \operatorname{tr}(M^T P_Z^\perp M)\}.$$

From (2.5) and (2.4), one has

$$(I-A)\mathbf{Y} = (I-Q_Z)(I - RE^+ R^T(I-Q_Z))\mathbf{Y}$$
$$= (P_Z^\perp + P_Z - Q_Z)(\boldsymbol{\eta} - \hat{\boldsymbol{\eta}} + Z\mathbf{b} + \boldsymbol{\varepsilon})$$
$$= P_Z^\perp(\boldsymbol{\eta} - \hat{\boldsymbol{\eta}}) + (P_Z - Q_Z)(\boldsymbol{\eta} - \hat{\boldsymbol{\eta}} + Z\mathbf{b} + \boldsymbol{\varepsilon}) + P_Z^\perp \boldsymbol{\varepsilon}$$
$$= P_Z^\perp(\boldsymbol{\eta} - \hat{\boldsymbol{\eta}}) + (P_Z - Q_Z)(\mathbf{Y} - \hat{\boldsymbol{\eta}}) + P_Z^\perp \boldsymbol{\varepsilon}.$$

It follows that

$$\mathbf{Y}^T(I-A)^2\mathbf{Y} = (\boldsymbol{\eta}-\hat{\boldsymbol{\eta}})^T P_Z^\perp(\boldsymbol{\eta}-\hat{\boldsymbol{\eta}}) + \boldsymbol{\varepsilon}^T P_Z^\perp \boldsymbol{\varepsilon} + 2(\boldsymbol{\eta}-\hat{\boldsymbol{\eta}})^T P_Z^\perp \boldsymbol{\varepsilon}$$
$$+ (\mathbf{Y}-\hat{\boldsymbol{\eta}})^T(P_Z-Q_Z)^2(\mathbf{Y}-\hat{\boldsymbol{\eta}}),$$

and hence

$$U(\lambda,\Sigma) - L_2(\lambda,\Sigma) - \frac{1}{n}\boldsymbol{\varepsilon}^T\boldsymbol{\varepsilon}$$
$$\text{(4.8)} \qquad = \frac{1}{n}(\mathbf{Y}-\hat{\boldsymbol{\eta}})^T(P_Z-Q_Z)^2(\mathbf{Y}-\hat{\boldsymbol{\eta}})$$
$$+ \frac{2}{n}(\boldsymbol{\eta}-\hat{\boldsymbol{\eta}})^T P_Z^\perp \boldsymbol{\varepsilon} - \frac{1}{n}\boldsymbol{\varepsilon}^T P_Z \boldsymbol{\varepsilon} + \frac{2}{n}\sigma^2 \operatorname{tr} A.$$

With an extra condition, $U(\lambda,\Sigma)$ and $V(\lambda,\Sigma)$ can be shown to track $L_2(\lambda,\Sigma)$ asymptotically.

CONDITION C.4. As $n \to \infty$, $R_1(\lambda,\Sigma) - R_2(\lambda,\Sigma) = o(R_1(\lambda,\Sigma))$.

Conditions C.2 and C.4 together imply that $R_1(\lambda,\Sigma) - R_2(\lambda,\Sigma) = o(R_2(\lambda,\Sigma))$ and $nR_2(\lambda,\Sigma) \to \infty$. Subtracting (4.7) from (4.3), some algebra yields

$$R_1(\lambda,\Sigma) - R_2(\lambda,\Sigma)$$
$$= \frac{1}{n}\boldsymbol{\eta}^T(I-M)^T(P_Z-Q_Z)^2(I-M)\boldsymbol{\eta}$$
$$\text{(4.9)} \qquad + \frac{1}{n}\operatorname{tr}\left(((P_Z-Q_Z) + (P_Z-Q_Z)RE^+R^T(P_Z-Q_Z))^2 ZBZ^T\right)$$
$$+ \frac{\sigma^2}{n}\operatorname{tr}\left((Q_Z + (P_Z-Q_Z)M)^T(Q_Z + (P_Z-Q_Z)M)\right).$$

The following lemma confirms the feasibility of Condition C.4 for fixed $p$.

LEMMA 4.3. *For fixed $p$, if (i) $\boldsymbol{\eta}^T(I - A(\lambda, \Sigma))P_Z(I - A(\lambda, \Sigma))\boldsymbol{\eta} = o(\boldsymbol{\eta}^T(I - A(\lambda, \Sigma))^2\boldsymbol{\eta})$, (ii) $\Sigma < \rho_n Z^T Z$ for $\rho_n^2 = o(R_1)$, and (iii) $\operatorname{tr}(Z^T Z)/n$ is bounded, then $R_1(\lambda, \Sigma) - R_2(\lambda, \Sigma) = o(R_1(\lambda, \Sigma))$.*

The proof of the lemma is given in Section 7. Condition (i) bars $(I - A)\boldsymbol{\eta}$ from being overloaded in the column space of $Z$; (ii) holds for $\Sigma$ of magnitude up to the order $O(\sqrt{n})$ when $Z^T Z$ grows at a rate $O(n)$, which is typical for fixed $p$. Alternatively, if $\rho_n = o(R_1)$ in (ii), which usually holds for bounded $\Sigma$, then (i) can be replaced by bounded $\boldsymbol{\eta}^T\boldsymbol{\eta}/n$; see the proof in Section 7.

THEOREM 4.3. *Under Conditions C.1, C.2 and C.4, as $n \to \infty$, one has*

$$U(\lambda, \Sigma) - L_2(\lambda, \Sigma) - \frac{1}{n}\boldsymbol{\varepsilon}^T\boldsymbol{\varepsilon} = o_p(L_2(\lambda, \Sigma)).$$

*If, in addition, Condition C.3 also holds, then*

$$V(\lambda, \Sigma) - L_2(\lambda, \Sigma) - \frac{1}{n}\boldsymbol{\varepsilon}^T\boldsymbol{\varepsilon} = o_p(L_2(\lambda, \Sigma)).$$

The proof of the theorem is given in Section 7.

Up to this point, we have considered purely real and purely latent random effects. In practice, one could have a mixture of real and latent random effects in the same setting. Partition $Z = (Z_1, Z_2)$ and $\mathbf{b}^T = (\mathbf{b}_1^T, \mathbf{b}_2^T)$ and assume $\mathbf{b}_1$ and $\mathbf{b}_2$ are independent so $B$ is block diagonal. Define

$$(4.10)\quad L_3(\lambda, \Sigma) = \frac{1}{n}(\hat{\boldsymbol{\eta}} + Z_1\hat{\mathbf{b}}_1 - \boldsymbol{\eta} - Z_1\mathbf{b}_1)^T P_{Z_2}^\perp (\hat{\boldsymbol{\eta}} + Z_1\hat{\mathbf{b}}_1 - \boldsymbol{\eta} - Z_1\mathbf{b}_1)$$

and $R_3(\lambda, \Sigma) = E[L_3(\lambda, \Sigma)]$, where $P_{Z_2}^\perp = I - Z_2(Z_2^T Z_2)^+ Z_2^T$; $L_3(\lambda, \Sigma)$ is a natural loss for $Z_1\mathbf{b}_1$ real and $Z_2\mathbf{b}_2$ latent. Replace $R_2(\lambda, \Sigma)$ in Condition C.4 by $R_3(\lambda, \Sigma)$.

CONDITION C.5. *As $n \to \infty$, $R_1(\lambda, \Sigma) - R_3(\lambda, \Sigma) = o(R_1(\lambda, \Sigma))$.*

A general result follows, of which the earlier theorems are special cases with nil $Z_1$ or nil $Z_2$.

THEOREM 4.4. *Under Conditions C.1, C.2 and C.5, as $n \to \infty$, one has*

$$U(\lambda, \Sigma) - L_3(\lambda, \Sigma) - \frac{1}{n}\boldsymbol{\varepsilon}^T\boldsymbol{\varepsilon} = o_p(L_3(\lambda, \Sigma)).$$

*If, in addition, Condition C.3 also holds, then*

$$V(\lambda, \Sigma) - L_3(\lambda, \Sigma) - \frac{1}{n}\boldsymbol{\varepsilon}^T\boldsymbol{\varepsilon} = o_p(L_3(\lambda, \Sigma)).$$

The proof of the theorem follows from straightforward modifications of the proof of Theorem 4.3 as given in Section 7.

**5. Empirical performance.** We now present simple simulations to illustrate the practical performance of generalized cross-validation in the context.

5.1. *Real random effects.* First consider a setting with real random effects covered by Theorems 4.1 and 4.2. One hundred replicates of samples were generated according to

$$Y_i = \eta(x_i) + b_{s_i} + \varepsilon_i, \qquad i = 1, \ldots, 100, \tag{5.1}$$

where $\eta(x) = 3\sin(2\pi x)$, $x_i$ a random sample from $U(0,1)$, $\varepsilon_i \sim N(0, 0.5^2)$, $b_s \sim N(0, 0.5^2)$ and $s_i \in \{1, \ldots, 10\}$, ten each. Cubic smoothing splines as described in Example 3.1 were calculated with $(\lambda_u, \Sigma_u)$ minimizing $U(\lambda, \Sigma)$ of (4.4), $(\lambda_v, \Sigma_v)$ minimizing $V(\lambda, \Sigma)$ of (4.1) and $(\lambda_m, \Sigma_m)$ minimizing $L_1(\lambda, \Sigma)$ of (4.2).

The loss $L_1(\lambda, \Sigma)$ was recorded for the fits. For the $V$ fit with $(\lambda_v, \Sigma_v)$, the variance estimate through

$$\hat{\sigma}^2 = \frac{\mathbf{Y}^T(I - A(\lambda_v, \Sigma_v))^2 \mathbf{Y}}{\text{tr}(I - A(\lambda_v, \Sigma_v))} \tag{5.2}$$

was also recorded; the variance estimate was proposed by Wahba [14] for independent data. The ratio $\sigma^2/\sigma_s^2$ as part of $\Sigma$ was "estimated" through $\Sigma_u$, $\Sigma_v$ or $\Sigma_m$.

It is known that cross-validation may lead to severe undersmoothing on up to about 10% replicates. To circumvent the problem, a simple modification proved to be very effective in the empirical studies of Kim and Gu [9]. The modified $V$ is given by

$$V_\alpha(\lambda, \Sigma) = \frac{n^{-1}\mathbf{Y}^T(I - A(\lambda, \Sigma))^2 \mathbf{Y}}{\{n^{-1}\text{tr}(I - \alpha A(\lambda, \Sigma))\}^2} \tag{5.3}$$

for some $\alpha > 1$. Similarly, $U$ can be modified by

$$U_\alpha(\lambda, \Sigma) = \frac{1}{n}\mathbf{Y}^T(I - A(\lambda, \Sigma))^2\mathbf{Y} + \frac{2}{n}\sigma^2 \alpha\, \text{tr}\, A(\lambda, \Sigma). \tag{5.4}$$

A good choice of $\alpha$ is around 1.4. The $U$ and $V$ fits with $\alpha = 1.2, 1.4, 1.6, 1.8$ were also calculated and the loss and variance estimates recorded.

The performances of $U_\alpha(\lambda, \Sigma)$ and $V_\alpha(\lambda, \Sigma)$ are illustrated in Figure 1. In the left and center frames, the losses $L_1(\lambda_u, \Sigma_u)$ and $L_1(\lambda_v, \Sigma_v)$ are plotted versus the minimum possible, for $\alpha = 1, 1.4$. The relative efficacy of $U_\alpha(\lambda, \Sigma)$ and $V_\alpha(\lambda, \Sigma)$ for $\alpha = 1, 1.2, 1.4, 1.6, 1.8$ is summarized in the right frame in box plots. Roughly speaking, $U_\alpha$ and $V_\alpha$ with $\alpha = 1$ are "unbiased" by Theorems 4.1 and 4.2, and setting $\alpha > 1$ introduces "bias." The top-tier performance may degrade slightly as $\alpha$ increases, but the worst cases are being pulled in for $\alpha$ up to $1.2 \sim 1.4$, where one appears to have the "minimax" performance.

Further details of the simulation are shown in Figure 2. In the left frame, $\lambda_u$ and $\lambda_v$ for $\alpha = 1$ and $\alpha = 1.4$ are plotted against each other, where a very small $\lambda$ by $\alpha = 1$ is seen to be pulled to the "normal" range by $\alpha = 1.4$. The number of cases with severe undersmoothing by cross-validation seems to be much less than what is typically seen in simulations with independent error; the phenomenon has yet to be understood. The center frame of Figure 2 plots the variance ratio $\sigma^2/\sigma_s^2$ "estimated" through $\Sigma_m$, $\Sigma_u$ and $\Sigma_v$. An interesting observation is the wide range of $\Sigma_m$, especially the many very small values, which effectively leave the term $\mathbf{z}^T\mathbf{b}$ unpenalized like the fixed effect terms in the null space of $J(\eta)$. The "estimates" through $\Sigma_u$ and $\Sigma_v$ appear far better in comparison, but remain highly unreliable. The upward trend of $\Sigma_u$ and $\Sigma_v$ with increasing $\alpha$ is somewhat expected, as larger $\alpha$ yields smoother estimates corresponding to larger penalty terms. In the right frame of Figure 2, the variance estimates by (5.2) are shown in box plots for $V$ fits with $\alpha = 1, 1.2, 1.4, 1.6, 1.8$, demonstrating generally adequate performance.

5.2. *Latent random effects.* For latent random effects, we keep the setting of (5.1) but replace $b_{s_i}$ by $b_{c_i}$, as in Example 3.3. One hundred replicates of samples were generated with $\eta(x_i)$ and $\varepsilon_i$ as in Section 5.1, and with $c_i \in \{1, 2\}$, 50 each, $b_1 \sim N(0, \sigma_1^2)$ for $\sigma_1^2 = 0.5^2$, and $b_2 \sim N(0, \sigma_2^2)$ for $\sigma_2^2 = 0.3^2$; the intracenter correlations are $0.25/(0.25 + 0.25) = 0.5$ for $c = 1$ and $0.09/(0.09 + 0.25) = 0.265$ for $c = 2$. Cubic smoothing splines were calculated with $\lambda$ and $\Sigma$ minimizing $U(\lambda, \Sigma)$, $V(\lambda, \Sigma)$ and $L_2(\lambda, \Sigma)$ of (4.6).

The simulation results are summarized in Figures 3 and 4. Figure 3 parallels Figure 1, except that $L_1(\lambda, \Sigma)$ is replaced by $L_2(\lambda, \Sigma)$. The left and center frames of Figure 4 summarize the "estimation" of the two parameters of $\Sigma$; note that the data contain only one "sample" from $N(0, \sigma_1^2)$ and one from $N(0, \sigma_2^2)$.

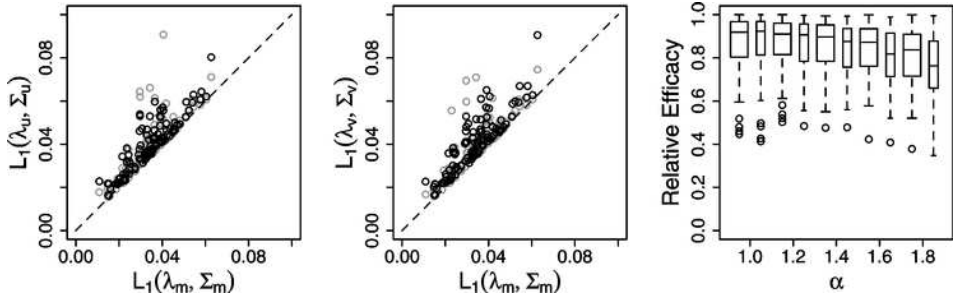

FIG. 1. *Simulation with real random effects.* Left *and* center: *Performances of* $U_\alpha(\lambda, \Sigma)$ *and* $V_\alpha(\lambda, \Sigma)$ *with* $\alpha = 1$ *(faded circles) and* $\alpha = 1.4$ *(circles).* Right: $L_1(\lambda_m, \Sigma_m)/L_1(\lambda_u, \Sigma_u)$ *(fatter boxes) and* $L_1(\lambda_m, \Sigma_m)/L_1(\lambda_v, \Sigma_v)$ *(thinner boxes) for* $\alpha = 1, 1.2, 1.4, 1.6, 1.8$.

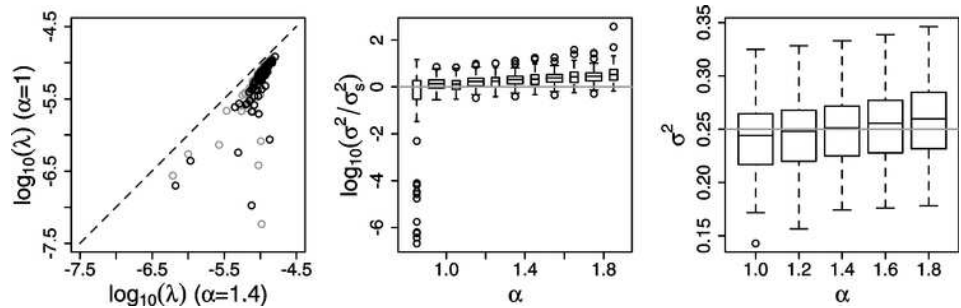

Fig. 2. *Simulation with real random effects.* Left: $\lambda_u$ *(faded circles) and* $\lambda_v$ *(circles) for* $\alpha = 1, 1.4$. *Center:* $\sigma^2/\sigma_s^2$ *"estimated" through* $\Sigma_m$ *(left thin box),* $\Sigma_u$ *(fatter boxes) and* $\Sigma_v$ *(thinner boxes).* Right: $\hat{\sigma}^2$. *The faded horizontal lines in center and right frames mark the true values.*

5.3. *Mixture random effects.* For mixture random effects, we simply add together $b_s$ of Section 5.1 and $b_c$ of Section 5.2, with the ten subjects nested under the two clusters, five each. One hundred replicates of samples were generated, with the specifications of $\eta(x)$, $\sigma^2$, $\sigma_s^2$, $\sigma_1^2$ and $\sigma_2^2$ remaining the same as in Sections 5.1 and 5.2. Cubic smoothing splines were calculated with $\lambda$ and $\Sigma$ minimizing $U(\lambda, \Sigma)$, $V(\lambda, \Sigma)$ and $L_3(\lambda, \Sigma)$ of (4.10). The counterpart of Figures 1 and 3 is shown in Figure 5. The "estimated" variance ratios are again highly unreliable, whereas $\hat{\sigma}^2$ demonstrates adequate performance, as seen in Figures 2 and 4; plots are omitted.

**6. Applications.** We now apply the technique to analyze a couple of real data sets.

6.1. *Tumor volume.* To study the sensitivity of a human prostate tumor to androgen deprivation, a preparation of the PC82 prostate cancer cell line was implanted under the skin of eight male nude mice. After 46 days, measurable tumors appeared on all eight mice; this day is referred to as day 0. On day 32, all mice were castrated. The tumors were measured roughly weekly over a 5-month period, resulting in 16 sets of measurements on the eight mice. Further details concerning the data can be found in [6], along with some analyses using parametric models.

We performed a nonparametric analysis of the data using the techniques developed. Taking the logarithm of the measured tumor volume as the response $Y$, the model of Example 3.1 was considered,

$$Y_i = \eta(x_i) + b_{s_i} + \varepsilon_i,$$

where $s = 1, \ldots, 8$. The exponential spline as discussed in Example 3.1 was used to estimate $\eta(x)$, but the generalized cross-validation score was minimized at $\theta = 0$, yielding a cubic spline fit. The fitted $\eta(x)$ is plotted in Figure 6 along with the data. The variance estimates are given by $\hat{\sigma}^2 = 0.1490$ and $\hat{\sigma}_s^2 = 0.0928$; remember that $\hat{\sigma}^2$ is trustworthy but $\hat{\sigma}_s^2$ can be grossly misleading, as shown in Section 5.

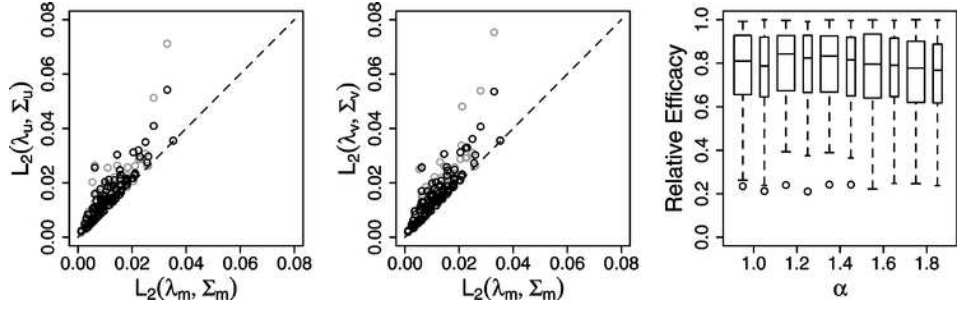

Fig. 3. *Simulation with latent random effects.* Left *and* center: *Performances of* $U_\alpha(\lambda, \Sigma)$ *and* $V_\alpha(\lambda, \Sigma)$ *with* $\alpha = 1$ *(faded circles) and* $\alpha = 1.4$ *(circles).* Right: $L_2(\lambda_m, \Sigma_m)/L_2(\lambda_u, \Sigma_u)$ *(fatter boxes) and* $L_2(\lambda_m, \Sigma_m)/L_2(\lambda_v, \Sigma_v)$ *(thinner boxes) for* $\alpha = 1, 1.2, 1.4, 1.6, 1.8$.

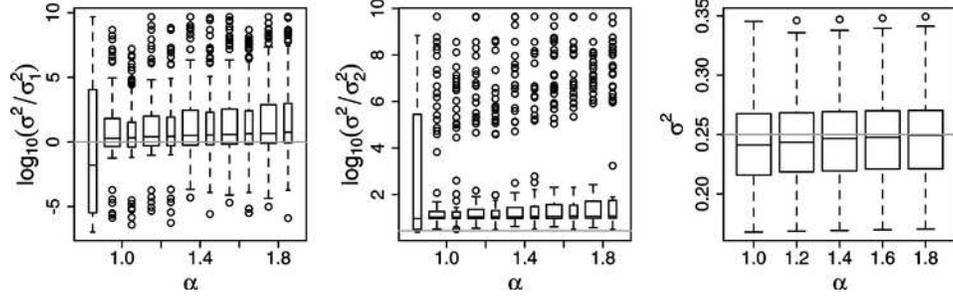

Fig. 4. *Simulation with latent random effects.* Left *and* center: $\sigma^2/\sigma_1^2$ *and* $\sigma^2/\sigma_2^2$ *"estimated" through* $\Sigma_m$ *(left thin box),* $\Sigma_u$ *(fatter boxes) and* $\Sigma_v$ *(thinner boxes).* Right: $\hat{\sigma}^2$. *The faded horizontal line marks the true values.*

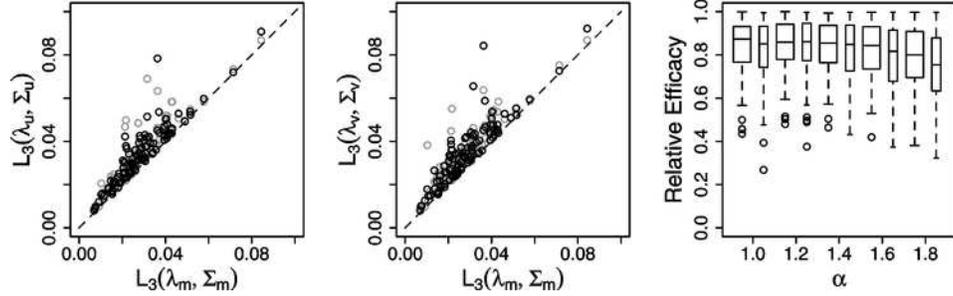

Fig. 5. *Simulation with mixture random effects.* Left *and* center: *Performances of* $U_\alpha(\lambda, \Sigma)$ *and* $V_\alpha(\lambda, \Sigma)$ *with* $\alpha = 1$ *(faded circles) and* $\alpha = 1.4$ *(circles).* Right: $L_3(\lambda_m, \Sigma_m)/L_3(\lambda_u, \Sigma_u)$ *(fatter boxes) and* $L_3(\lambda_m, \Sigma_m)/L_3(\lambda_v, \Sigma_v)$ *(thinner boxes) for* $\alpha = 1, 1.2, 1.4, 1.6, 1.8$.

6.2. *Treatment of multiple sclerosis.* A randomized, double-blind clinical trial was conducted to study the treatment of multiple sclerosis by azathioprine (AZ) and methylprednisolone (MP). Patients were assigned randomly to three groups: (i) the PP group receiving placebos for both AZ and MP, (ii) the AP group receiving real AZ and placebo MP; and (iii) the AM group receiving real AZ and MP. The abundance of lymphocytes bearing a protein called $F_C$ receptor was measured in the form of the so-called AFCR levels. Blood samples were drawn prior to the initiation of therapy, at the initiation, in weeks 4, 8 and 12, and every 12 weeks thereafter for the remainder of the trial. A total of 48 patients were represented in the data, with 17 on PP, 15 on AP and 16 on AM. There were "missing" values in the sense that blood samples were not drawn from all patients at every time point. Detailed descriptions of the study can be found in [7] and further references therein. A analysis of the data using parametric models was conducted by Heitjan [7].

We now present a nonparametric analysis of the data using the formulation of Example 3.2. Following [7], the responses $Y_i$ are taken as the square roots of the AFCR measures. The model is of the form

$$Y_i = \eta(x_i, \tau_i) + b_{s_i} + \varepsilon_i,$$

where the patient identification $s$ is nested under the treatment level $\tau$. The "missing" values pose no problem for our treatment. The fitted cubic splines are plotted in Figure 7 with the data superimposed. The smoothing parameter $\theta_{1,2}$ was effectively set to 0 by cross-validation, so the interaction $\eta_{1,2}(x,\tau)$ consists of only parametric terms with the basis $(I_{[\tau=j]} - 1/3)x$, $j = 1, 2$; see Example 3.2 for the notation. The variance estimates were given by $\hat{\sigma}^2 = 12.81$ and $\hat{\sigma}_s^2 = 6.624$.

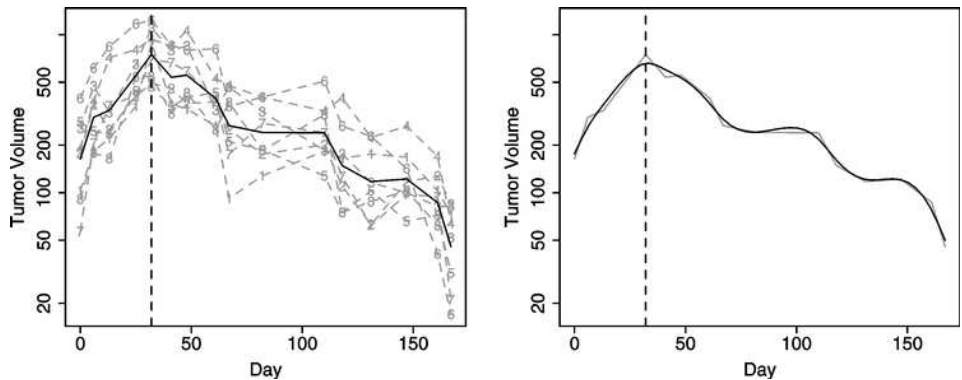

Fig. 6. *Cubic spline fits of tumor volume.* Left: *Tumor volume measurements (dashed lines) and their geometric mean (solid line).* Right: *Fitted $\eta(x)$ (solid line), with the geometric mean of measurements superimposed (faded line). The castration time is marked by the vertical line.*

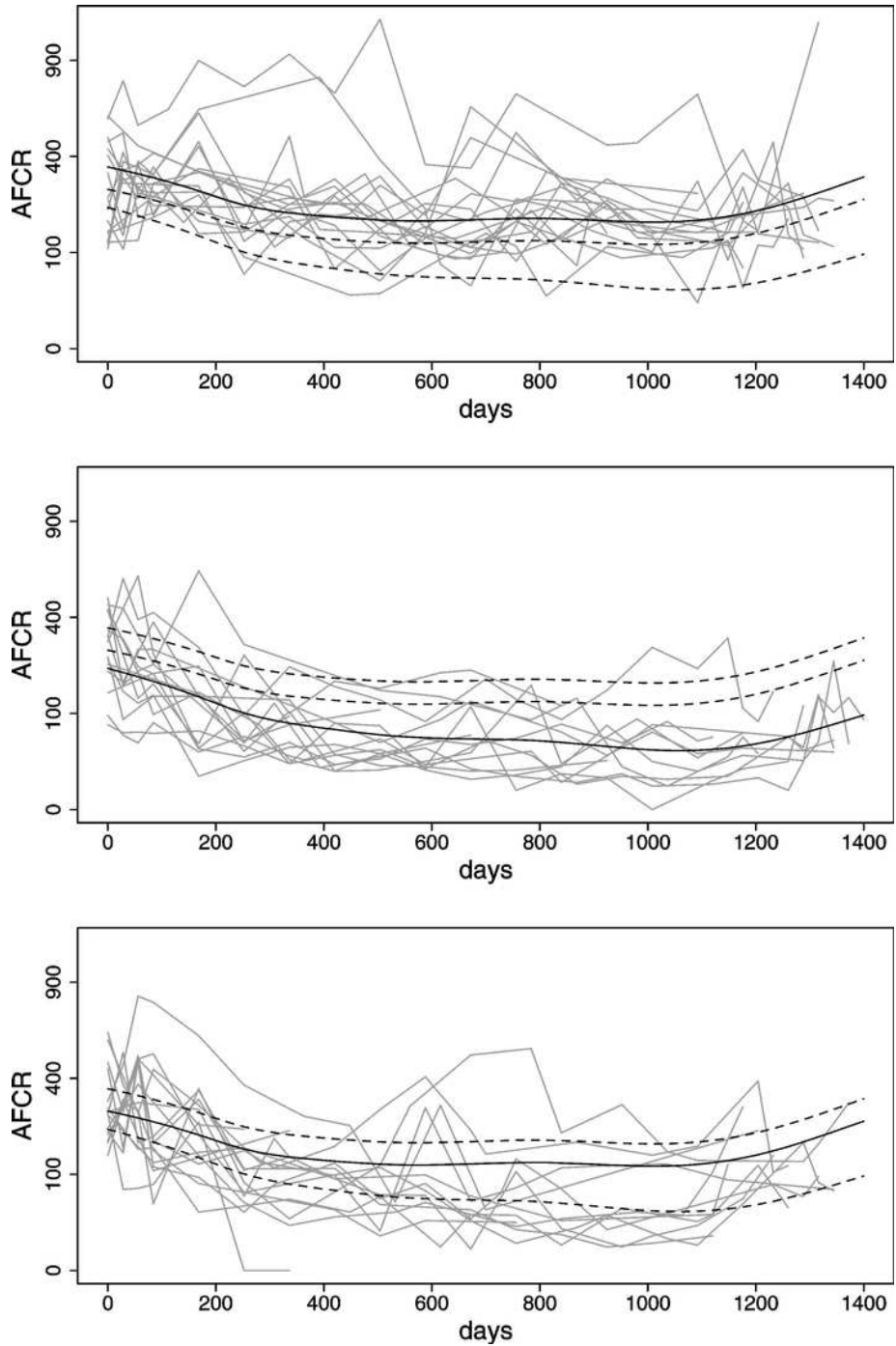

FIG. 7. *Cubic spline fits of AFCR levels. From top to bottom: the PP, AP and AM groups. The fitted $\eta(x,\tau)$, $\tau = PP, AP, AM$, are in solid lines in their respective frames, with the corresponding data superimposed as faded lines and the other two estimates as dashed lines.*

**7. Proofs.** This section collects the proofs of the lemmas and theorems of Section 4. The following lemmas govern some of the calculations.

LEMMA 7.1. *For $Z\mathbf{b} \sim N(\mathbf{0}, ZBZ^T)$ and $\boldsymbol{\varepsilon} \sim N(\mathbf{0}, \sigma^2 I)$, independent of each other, one has*

$$\mathrm{Var}[\mathbf{b}^T Z^T C Z \mathbf{b}] = 2\,\mathrm{tr}(CZBZ^T C^T ZBZ^T),$$
$$\mathrm{Var}[\boldsymbol{\varepsilon}^T C \boldsymbol{\varepsilon}] = 2\sigma^4 \,\mathrm{tr}(CC^T),$$
$$\mathrm{Var}[\mathbf{b}^T Z^T C \boldsymbol{\varepsilon}] = \sigma^2 \,\mathrm{tr}(CC^T ZBZ^T).$$

The proof of the lemma is straightforward.

LEMMA 7.2. *For $M = RE^+ R^T (I - Q_Z)$, where $E = R^T (I - Q_Z) R + n\lambda Q$, one has*

$$M^T P_Z^\perp M + M^T (P_Z - Q_Z) M$$
$$= M^T (I - Q_Z) M \leq I,$$
$$(I - M)^T P_Z^\perp (I - M) + (I - M)^T (P_Z - Q_Z)(I - M)$$
$$= (I - M)^T (I - Q_Z)(I - M) \leq 4I.$$

PROOF. It is straightforward to show that $M^T (I - Q_Z) M \leq I$. Now for an arbitrary vector $\mathbf{x}$,

$$\mathbf{x}^T (I - M)^T (I - Q_Z)(I - M)\mathbf{x}$$
$$= \mathbf{x}^T (I - Q_Z)\mathbf{x} + \mathbf{x}^T M^T (I - Q_Z) M \mathbf{x} - 2\mathbf{x}^T (I - Q_Z) M \mathbf{x}^T \leq 4\mathbf{x}^T \mathbf{x},$$

where the Cauchy–Schwarz inequality is used to bound the cross term. □

Also note that $B$ is fixed, thus having bounded eigenvalues, and that $X^T X$ and $XX^T$ share nonzero eigenvalues for all matrices $X$.

We are now ready for the proofs of the lemmas and theorems of Section 4.

PROOF OF LEMMA 4.1. Recall from (4.3),

$$R_1(\lambda, \Sigma) = \frac{1}{n}\boldsymbol{\eta}^T (I - A)^2 \boldsymbol{\eta} + \frac{1}{n}\,\mathrm{tr}((I - A)^2 ZBZ^T) + \frac{\sigma^2}{n}\,\mathrm{tr}\, A^2.$$

Using (2.7), the first term is seen to be of the order $O(\lambda^s)$, and the third term is of order $O(n^{-1}\lambda^{-1/r} + n^{-1}p)$. Again by (2.7),

$$(I - A)Z = (I - \tilde{A})Z(I - (Z^T (I - \tilde{A})Z + \Sigma)^{-1} Z^T (I - \tilde{A})Z)$$
$$= (I - \tilde{A})Z(Z^T (I - \tilde{A})Z + \Sigma)^{-1}\Sigma,$$

thus
$$Z^T(I-A)^2 Z \le \Sigma(Z^T(I-\tilde{A})Z + \Sigma)^{-1}\Sigma,$$

so Condition C.1 implies an upper bound on the eigenvalues of $(I-A)ZBZ^T(I-A)$, and the second term is of order $O(n^{-1}p)$. The proof is complete. $\square$

PROOF OF THEOREM 4.1. In light of (4.5), it suffices to show that

(7.1) $$L_1(\lambda,\Sigma) - R_1(\lambda,\Sigma) = o_p(R_1(\lambda,\Sigma)),$$

(7.2) $$n^{-1}(\boldsymbol{\eta} + Z\mathbf{b})^T(I-A)\boldsymbol{\varepsilon} = o_p(R_1(\lambda,\Sigma)),$$

(7.3) $$n^{-1}(\boldsymbol{\varepsilon}^T A \boldsymbol{\varepsilon} - \sigma^2 \operatorname{tr} A) = o_p(R_1(\lambda,\Sigma)).$$

To see (7.1), note that

$$\operatorname{Var}[L_1(\lambda,\Sigma)] = n^{-2}\operatorname{Var}[2\boldsymbol{\eta}^T(I-A)^2 Z\mathbf{b} - 2\boldsymbol{\eta}^T(I-A)A\boldsymbol{\varepsilon} \\ + \mathbf{b}^T Z^T(I-A)^2 Z\mathbf{b} - 2\mathbf{b}^T Z^T(I-A)A\boldsymbol{\varepsilon} + \boldsymbol{\varepsilon}^T A^2\boldsymbol{\varepsilon}].$$

Since Condition C.1 implies an upper bound on the eigenvalues of $(I-A)ZBZ^T(I-A)$, one has

$$n^{-2}\operatorname{Var}[\boldsymbol{\eta}^T(I-A)^2 Z\mathbf{b}] = n^{-2}\boldsymbol{\eta}^T(I-A)^2 ZBZ^T(I-A)^2\boldsymbol{\eta} \\ = n^{-1}O(R_1) = o(R_1^2),$$

where the last equation is by Condition C.2. Likewise,

$$n^{-2}\operatorname{Var}[\boldsymbol{\eta}^T(I-A)A\boldsymbol{\varepsilon}] = n^{-2}\sigma^2\boldsymbol{\eta}^T(I-A)A^2(I-A)\boldsymbol{\eta} = o(R_1^2),$$
$$n^{-2}\operatorname{Var}[\mathbf{b}^T Z^T(I-A)^2 Z\mathbf{b}] = 2n^{-2}\operatorname{tr}((I-A)^2 ZBZ^T(I-A)^2 ZBZ^T) = o(R_1^2),$$
$$n^{-2}\operatorname{Var}[\mathbf{b}^T Z^T(I-A)A\boldsymbol{\varepsilon}] = n^{-2}\sigma^2\operatorname{tr}((I-A)A^2(I-A)ZBZ^T) = o(R_1^2),$$
$$n^{-2}\operatorname{Var}[\boldsymbol{\varepsilon}A^2\boldsymbol{\varepsilon}] = 2n^{-2}\sigma^4\operatorname{tr}A^4 = o(R_1^2).$$

Summing up, and bounding the covariances between the terms by the Cauchy–Schwarz inequality, one has $\operatorname{Var}[L_1(\lambda,\Sigma)] = o(R_1^2(\lambda,\Sigma))$, and hence (7.1). Similar calculations yield (7.2) and (7.3), completing the proof. $\square$

PROOF OF LEMMA 4.2. From (2.7), one has $\operatorname{tr} A \le \operatorname{tr}\tilde{A} + p$ and $\operatorname{tr} A^2 \ge \operatorname{tr}\tilde{A}^2$, so

$$\frac{(n^{-1}\operatorname{tr} A)^2}{n^{-1}\operatorname{tr} A^2} \le \frac{(n^{-1}\operatorname{tr}\tilde{A} + n^{-1}p)^2}{n^{-1}\operatorname{tr}\tilde{A}^2} = O(n^{-1}\lambda^{-1/r} + n^{-1}p + n^{-1}p^2\lambda^{1/r}).$$

The lemma follows as $\lambda \to 0$ and $n\lambda^{1/r} \to \infty$. $\square$

PROOF OF LEMMA 4.3. Recall from (4.9) that

$$R_1(\lambda, \Sigma) - R_2(\lambda, \Sigma)$$
$$= \frac{1}{n}\boldsymbol{\eta}^T(I-M)^T(P_Z - Q_Z)^2(I-M)\boldsymbol{\eta}$$
$$+ \frac{1}{n}\text{tr}(((P_Z - Q_Z) + (P_Z - Q_Z)RE^+R^T(P_Z - Q_Z))^2 ZBZ^T)$$
$$+ \frac{\sigma^2}{n}\text{tr}((Q_Z + (P_Z - Q_Z)M)^T(Q_Z + (P_Z - Q_Z)M)).$$

Since $D = Z^TZ + \Sigma < (1+\rho_n)Z^TZ$, one has $P_Z - Q_Z < \rho_n P_Z/(1+\rho_n) < \rho_n P_Z$. For the first line, noting that $P_Z - Q_Z = P_Z(I - Q_Z)$ and $(I-Q_Z)(I-M) = I - A$, one has

$$\frac{1}{n}\boldsymbol{\eta}^T(I-M)^T(P_Z - Q_Z)^2(I-M)\boldsymbol{\eta} = \frac{1}{n}\boldsymbol{\eta}^T(I-A)P_Z(I-A)\boldsymbol{\eta} = o(R_1).$$

Alternatively, with $\boldsymbol{\eta}^T\boldsymbol{\eta}/n$ bounded,

$$\frac{1}{n}\boldsymbol{\eta}^T(I-M)^T(P_Z - Q_Z)^2(I-M)\boldsymbol{\eta}$$
$$\leq \rho_n \frac{1}{n}\boldsymbol{\eta}^T(I-M)^T(P_Z - Q_Z)(I-M)\boldsymbol{\eta}$$
$$= O(\rho_n)$$

as $(I-M)^T(P_Z - Q_Z)(I-M) \leq 4I$. For the second line, note that

$$\frac{1}{n}\text{tr}((P_Z - Q_Z)ZBZ^T(P_Z - Q_Z)) \leq \rho_n^2 \frac{1}{n}\text{tr}(ZBZ^T) = o(R_1)$$

and, with $F = (P_Z - Q_Z)^{1/2}RE^+R^T(P_Z - Q_Z)^{1/2} \leq I$ and hence $F(P_Z - Q_Z)F \leq \rho_n I$, that

$$\frac{1}{n}\text{tr}(((P_Z - Q_Z)RE^+R^T(P_Z - Q_Z))^2 ZBZ^T)$$
$$= \frac{1}{n}\text{tr}(B^{1/2}Z^T(P_Z - Q_Z)^{1/2}F(P_Z - Q_Z)F(P_Z - Q_Z)^{1/2}ZB^{1/2})$$
$$\leq \rho_n \frac{1}{n}\text{tr}(B^{1/2}Z^T(P_Z - Q_Z)ZB^{1/2}) = \rho_n^2 \frac{1}{n}\text{tr}(ZBZ^T) = o(R_1);$$

the cross term can be bounded by the Cauchy–Schwarz inequality. For the third line, note that $n^{-1}\text{tr}\,Q_Z^2 \leq p/n = o(R_1)$ and that $M^T(P_Z - Q_Z)^2 M \leq I$ has no more than $p$ nonzero eigenvalues. The proof is now complete. □

PROOF OF THEOREM 4.3. Recall from (4.7) that

$$R_2(\lambda, \Sigma) = \frac{1}{n}\{\boldsymbol{\eta}^T(I-M)^T P_Z^\perp(I-M)\boldsymbol{\eta}$$
$$+ \text{tr}(M^T P_Z^\perp M Z B Z^T) + \sigma^2 \text{tr}(M^T P_Z^\perp M)\}.$$

Plugging $\hat{\boldsymbol{\eta}} = M(\boldsymbol{\eta} + Z\mathbf{b} + \boldsymbol{\varepsilon})$ into (4.8) and grouping terms, some algebra leads to

$$U(\lambda, \Sigma) - L_2(\lambda, \Sigma) - \frac{1}{n}\boldsymbol{\varepsilon}^T\boldsymbol{\varepsilon}$$

$$= \frac{1}{n}(\boldsymbol{\eta} + Z\mathbf{b})^T(I - M)^T(P_Z - Q_Z)^2(I - M)(\boldsymbol{\eta} + Z\mathbf{b})$$

$$+ \frac{2}{n}\boldsymbol{\eta}^T(I - M)^T(P_Z - Q_Z)^2(I - M)\boldsymbol{\varepsilon} + \frac{2}{n}\boldsymbol{\eta}^T(I - M)^T P_Z^\perp \boldsymbol{\varepsilon}$$

(7.4)

$$+ \frac{2}{n}\mathbf{b}^T Z^T(I - M)^T(P_Z - Q_Z)^2(I - M)\boldsymbol{\varepsilon} - \frac{2}{n}\mathbf{b}^T Z^T M^T P_Z^\perp \boldsymbol{\varepsilon}$$

$$+ \frac{1}{n}\boldsymbol{\varepsilon}^T(Q_Z + (P_Z - Q_Z)M)^T(Q_Z + (P_Z - Q_Z)M)\boldsymbol{\varepsilon}$$

$$- \frac{1}{n}(\boldsymbol{\varepsilon}^T A \boldsymbol{\varepsilon} - \sigma^2 \operatorname{tr} A).$$

To prove the first part of the theorem, it suffices to show that (7.4) is of order $o_p(R_2(\lambda, \Sigma))$ and that

(7.5) $\qquad L_2(\lambda, \Sigma) - R_2(\lambda, \Sigma) = o_p(R_2(\lambda, \Sigma)).$

Taking the expectation of the first line of (7.4), one has

$$\frac{1}{n}E[(\boldsymbol{\eta} + Z\mathbf{b})^T(I - M)^T(P_Z - Q_Z)^2(I - M)(\boldsymbol{\eta} + Z\mathbf{b})]$$

$$= \frac{1}{n}\boldsymbol{\eta}(I - M)^T(P_Z - Q_Z)^2(I - M)\boldsymbol{\eta}$$

$$+ \frac{1}{n}\operatorname{tr}(((P_Z - Q_Z) - (P_Z - Q_Z)RE^+R^T(P_Z - Q_Z))^2 ZBZ^T)$$

$$= O(R_1 - R_2) = o(R_2),$$

where Condition C.4 is used. Similarly, the expectation of the fourth line of (7.4) gives

$$\frac{1}{n}E[\boldsymbol{\varepsilon}^T(Q_Z + (P_Z - Q_Z)M)^T(Q_Z + (P_Z - Q_Z)M)\boldsymbol{\varepsilon}]$$

$$= \frac{\sigma^2}{n}\operatorname{tr}((Q_Z + (P_Z - Q_Z)M)^T(Q_Z + (P_Z - Q_Z)M))$$

$$= O(R_1 - R_2) = o(R_2).$$

For the two terms on the second line of (7.4), noting that $(I - M)^T(P_Z - Q_Z)^2(I - M) \leq 4I$,

$$n^{-2}\operatorname{Var}[\boldsymbol{\eta}^T(I - M)^T(P_Z - Q_Z)^2(I - M)\boldsymbol{\varepsilon}]$$

$$\leq 4n^{-2}\sigma^2 \boldsymbol{\eta}^T(I - M)^T(P_Z - Q_Z)^2(I - M)\boldsymbol{\eta} = o(R_2^2)$$

by Conditions C.2 and C.4, and

$$n^{-2} \operatorname{Var}[\boldsymbol{\eta}^T (I - M)^T P_Z^\perp \boldsymbol{\varepsilon}] = n^{-2} \sigma^2 \boldsymbol{\eta}^T (I - M)^T P_Z^\perp (I - M) \boldsymbol{\eta}$$
$$= n^{-1} O(R_2) = o(R_2^2).$$

Likewise, the third-line terms in (7.4) give

$$n^{-2} \operatorname{Var}[\mathbf{b}^T Z^T (I - M)^T (P_Z - Q_Z)^2 (I - M) \boldsymbol{\varepsilon}]$$
$$\leq 2 n^{-2} \sigma^2 \operatorname{tr}(((P_Z - Q_Z) - (P_Z - Q_Z) R E^+ R^T (P_Z - Q_Z))^2 Z B Z^T)$$
$$= o(R_2^2),$$

and

$$n^{-2} \operatorname{Var}[\mathbf{b}^T Z^T M^T P_Z^\perp \boldsymbol{\varepsilon}] = 2 n^{-2} \sigma^2 \operatorname{tr}(M^T P_Z^\perp M Z B Z^T) = n^{-1} O(R_2) = o(R_2^2).$$

The fifth line of (7.4) is (7.3), which is of order $o_p(R_1) = o_p(R_2)$ by Condition C.4. To see (7.5), note that

$$\operatorname{Var}[L_2(\lambda, \Sigma)]$$
$$= n^{-2} \operatorname{Var}[2 \boldsymbol{\eta}^T (M - I)^T P_Z^\perp M Z \mathbf{b} + 2 \boldsymbol{\eta}^T (M - I)^T P_Z^\perp M \boldsymbol{\varepsilon}$$
$$+ \mathbf{b}^T Z^T M^T P_Z^\perp M Z \mathbf{b} + \mathbf{b}^T Z^T M^T P_Z^\perp M \boldsymbol{\varepsilon} + \boldsymbol{\varepsilon}^T M^T P_Z^\perp M \boldsymbol{\varepsilon}].$$

Using (2.6), one has $MZ = \tilde{A} Z (Z^T (I - \tilde{A}) Z + \Sigma)^{-1} \Sigma$, $P_Z^\perp MZ = P_Z^\perp (\tilde{A} - I) Z (Z^T (I - \tilde{A}) Z + \Sigma)^{-1} \Sigma$, thus $Z^T M^T P_Z^\perp MZ \leq \Sigma (Z^T (I - \tilde{A}) Z + \Sigma)^{-1} \Sigma$, so Condition C.1 implies bounded eigenvalues for $P_Z^\perp MZBZ^T M^T P_Z^\perp$. It then follows that

$$n^{-2} \operatorname{Var}[\boldsymbol{\eta}^T (M - I)^T P_Z^\perp M Z \mathbf{b}] = n^{-2} \boldsymbol{\eta}^T (I - M)^T P_Z^\perp M Z B Z^T M^T P_Z^\perp (I - M) \boldsymbol{\eta}$$
$$= o(R_2^2),$$
$$n^{-2} \operatorname{Var}[\boldsymbol{\eta}^T (M - I)^T P_Z^\perp M \boldsymbol{\varepsilon}] = n^{-2} \sigma^2 \boldsymbol{\eta}^T (I - M)^T P_Z^\perp M M^T P_Z^\perp (I - M) \boldsymbol{\eta}$$
$$= o(R_2^2),$$
$$n^{-2} \operatorname{Var}[\mathbf{b}^T Z^T M^T P_Z^\perp M Z \mathbf{b}] = 2 n^{-2} \operatorname{tr}(M^T P_Z^\perp M Z B Z^T M^T P_Z^\perp M Z B Z^T)$$
$$= o(R_2^2),$$
$$n^{-2} \operatorname{Var}[\mathbf{b}^T Z^T M^T P_Z^\perp M \boldsymbol{\varepsilon}] = n^{-2} \sigma^2 \operatorname{tr}(M^T P_Z^\perp M Z B Z^T M^T P_Z^\perp M) = o(R_2^2),$$
$$n^{-2} \operatorname{Var}[\boldsymbol{\varepsilon}^T M^T P_Z^\perp M \boldsymbol{\varepsilon}] = 2 n^{-2} \sigma^4 \operatorname{tr}(M^T P_Z^\perp M M^T P_Z^\perp M) = o(R_2^2).$$

Collecting terms and bounding the covariances between the terms by the Cauchy–Schwarz inequality, one has $\operatorname{Var}[L_2(\lambda, \Sigma)] = o(R_2^2(\lambda, \Sigma))$, and hence (7.5). The proof of the first part of the theorem is now complete.

Given the first part of the theorem, the second part follows from the proof of Theorem 3.3 in [2], page 66. □

**8. Discussion.** In this article we studied the optimal smoothing of nonparametric mixed-effect models through generalized cross-validation. The asymptotic analysis was backed by simulation studies with sample size as small as 100. Related practical issues such as variance estimation were also explored in the simulation studies. As a sequel to this work, the optimal smoothing of non-Gaussian longitudinal data has been studied in [4] on an empirical basis. The methods have been implemented in the open-source R package gss by the first author.

While many correlated errors can be cast as variance components with low-rank random effects, some others do not conform, which spells the limitation of the techniques developed here; an important nonconforming case is serial or spatial correlation. On the flip side, the nonparametric $\eta(x)$ can be interpreted as a realization of a Gaussian process under the Bayes model of a smoothing spline, so there remains a potential identifiability problem of some sort between $\eta(x)$ and a separate serial or spatial correlation, unless the serial or spatial correlation is independent of $x$. Optimal smoothing for penalized likelihood estimation with serially or spatially correlated data is treated in a recent study by Gu and Han [3].

**Acknowledgments.** The authors are indebted to two referees, whose comments and suggestions led to a much improved presentation. The authors are also grateful to a co-editor and a third referee, who spotted a sloppy but critical technical flaw in an earlier version of the article.

Department of Statistics  
Purdue University  
West Lafayette, Indiana 47907  
USA  
e-mail: chong@stat.purdue.edu

Department of Statistics  
Harvard University  
Cambridge, Massachusetts 02138  
USA  
e-mail: pingma@stat.harvard.edu